\documentclass[11pt]{amsart}
\numberwithin{equation}{section}

\usepackage{amsmath,amsthm,amssymb,amsfonts,array,latexsym,amsthm,pdfsync,hyperref, enumerate, ulem, bbm}
\usepackage[usenames,dvipsnames]{color}
\theoremstyle{plain}

\newtheorem{thm}{Theorem}[section]

\newtheorem{cor}[thm]{Corollary}
\newtheorem{claim}[thm]{Claim}
\newtheorem{conj}[thm]{Conjecture}

\newtheorem{defin}[thm]{Definition}
\newtheorem{rem}[thm]{Remark}

\theoremstyle{remark}

\newcommand{\RR}{\ensuremath{\mathbb R}}

\newcommand{\calT}{\ensuremath{\mathcal T}}
\newcommand{\calW}{\ensuremath{\mathcal W}}

\date{\today}

\title{An Elementary Proof of the 3 Dimensional Simplex Mean Width Conjecture}
\author{Aaron Goldsmith}
\address{Aaron Goldsmith\\ Unaffiliated}
\email{aaron.goldsmith163@gmail.com}

\newcommand{\sphd}{\mathbb{S}^{d-1}}
\newcommand{\sph}{\mathbb{S}^2}

\begin{document}

\pagestyle{plain}
\maketitle

\begin{abstract}
After a Hessian computation, we quickly prove the 3D simplex mean width conjecture using classical methods in section \ref{dim3}. Then, we generalize some components to $d$ dimensions in section \ref{dimd}.
\end{abstract}

\section{Introduction}
Let $B_2^d\subset \RR^d$ be the standard Euclidean ball, with basis $(e_1,\dots,e_d)$. Denote $\sphd:=\partial B_2^d$. Let $\delta(X,Y)$ be the arclength on $\sphd$ and $\mu$ be the uniform probability measure on $\sphd$.\\

The support function of a convex body $K\subset B_2^d$ is

\[h_K(u):=\max_{x\in K}u\cdot x\]

and the mean width

\[w(K):=2\int_{\sphd} h_K(u)d\mu(u)\]

{\conj\label{SMWC} (Simplex Mean Width Conjecture) Of all simplices contained in $B_2^d$, the inscribed regular simplex has the maximum mean width, and is unique up to isometry.}\\

Conjecture \ref{SMWC} was mentioned in the survey by Gritzmann and Klee \cite{GK} and by Klee several times in his talks. Litvak surveyed the problem more recently in arXiv:math/0606350 [math.DG]. It's related to the problem of recovering transmissions from a noisy signal and has been assumed to be true by information theorists \cite{AVB1, AVB2}. The $d+1$ dimensional Gaussian Random Vector Maximum Conjecture is equivalent and was proved in 4 dimensions arXiv:2008.04827v2 [math.PR]. In section \ref{dim3}, we give a simpler and more classical proof. The following claim will be useful throughout.

{\claim\label{basic} In order that $\triangle:=\mathrm{conv}(v_0,\dots, v_d)$ maximizes $w(\cdot)$ over all simplices contained in $B_2^d$, it must be that 

\begin{enumerate}[(a)]
\item $B_2^d$ is the smallest ball containing $\triangle$.
\item $v_i\in\sphd$ for $0\le i\le n$, i.e. $\sphd$ is the circumsphere of $\triangle$.
\item The closed hemispheres centered at $v_i$ cover $\sphd$.
\end{enumerate}}

\begin{proof}
 \begin{enumerate}[(a)]
\item If $\triangle\subset Q+rB_2^d$ with $r<1$, then $(\triangle-Q)/r$ is a simplex contained in $B_2^d$ and

\[w((\triangle-Q)/r)=w(\triangle/r)=w(\triangle)/r>w(\triangle)\]

\item Suppose that $|v_0|<1$. Extend the edge $v_1v_0$ past $v_0$ until it meets $\sphd$ at a point $v_{n+1}$. The simplex $\tilde{\triangle}:=\mathrm{conv}\{v_1,\dots,v_{n+1}\}\subset B_2^d$ contains $v_0$ as a convex combination of $v_{n+1}$ and $v_1$, so it strictly contains $\triangle$. Since $w$ is strictly increasing under set inclusion, $w(\triangle)<w(\tilde{\triangle})$. Therefore, if $|v_i|<1$ for any $i$, then $w(\triangle)$ is not maximum. 
\item Suppose $v\in\sphd$ is such that $\max v\cdot v_i =m<0$ for all $i$. It follows that each $v_i$ is contained in the spherical cap $\{v\cdot x\le m\}$, which has diameter less than $1$. Part (a) implies there is no such $v$. 
\end{enumerate}
\end{proof}

The support function of the simplex $\triangle=\mathrm{conv}\{v_0,\dots, v_d\}$ reduces to a maximum over the vertex set:

\[\max_{x\in \triangle}u\cdot x=\max_{0\le i\le d} u\cdot v_i\]

As such, for each $i$ define the $i$th Voronoi cell on the sphere to be

\[S_i:=\{u\in\sphd:u\cdot v_i= \max_{x\in\triangle}u\cdot x\}\]

The mean width of a simplex is then

\begin{align}w(\triangle)=&2\int_{\sphd}\max_{0\le i\le d}X\cdot v_i d\mu(X)\nonumber\\ =& 2\sum_{i=0}^d \int_{S_i}X\cdot v_i d\mu(X)\label{SumMAR}\end{align}

Now define
\begin{align}M_A(R):=\int_R X\cdot A\ d\mu(X)\label{IntROnS}\end{align}

\section{SMWC dimension 3}\label{dim3}

Capital Roman letters $A,B,C$ represent either vertices of a spherical triangle or their angles, depending on context; small Roman letters $a,b,c,$ are used for either the opposite legs or their arclength. Next, we present some classical results about trigonometry in $\mathbb{S}^2$ (see \cite{IT}).

{\thm[Spherical Law of Sines]\label{SLOS} For a general triangle $\triangle ABC\subset \sph$,
\begin{align*}\frac{\sin a}{\sin A}=\frac{\sin b}{\sin B}=\frac{\sin c}{\sin C}\end{align*}}

{\thm[Napier] Suppose $\triangle ABC$ is a right spherical triangle. Then,	
\begin{eqnarray}
\cos c=&\cos a\cos b\label{SPT}\\
\cos A=&\cos a\sin B\label{Nap1}\\
\cos A=&\tan b\cot c\label{Nap2}\\
\tan a=&\tan A\sin b\label{Nap3}
\end{eqnarray} and, of course, the same if $(A,a)$ and $(B,b)$ are interchanged.}\\

Claim \ref{parama} gives a parametrization of an arc as measured from the opposite vertex and Theorem \ref{avgheight3D} a simple formula reminiscent of the area of a Euclidean triangle.

{\claim\label{parama} Let $\triangle ABC\subset \sph$ be a right spherical triangle with $A=(1,0,0)$ and the right angled vertex $C=(\cos b,\sin b, 0)$. Side $a$ is parametrized by $(\cos\Phi(\theta),\cos\theta\sin\Phi(\theta),\sin\theta\sin\Phi(\theta))$  as $\theta\in[0, A]$ where \[\sin^2\Phi(\theta)=\frac{\tan^2 b}{\tan^2 b+\cos^2\theta}\]}

\begin{proof}Napier's rule \ref{Nap2} says $\cos\theta=\tan b\cot \Phi(\theta)$, so

\begin{align*}\sin^2\Phi(\theta)=\frac{1}{1+\cot^2\Phi(\theta)}=\frac{\tan^2 b}{\tan^2 b+\cos^2\theta}\end{align*}
\end{proof}

{\thm\label{avgheight3D} Let $\triangle ABC$ be a right spherical triangle with $C=\pi/2$. Then,

\begin{align*}4\pi\sigma(\triangle ABC)M_A(\triangle ABC)=\frac{1}{2}a\sin b\end{align*}}

\begin{proof}Set $A=(1,0,0)$. Uniform measure on $\sph$ is $d\mu=(4\pi)^{-1}\sin\varphi d\varphi d\theta$, giving 

\begin{align*}
8\pi\int_{\triangle ABC} X\cdot A d\mu(X)=&\int_0^A\int_0^{\Phi(\theta)}\cos\varphi\sin\varphi d\varphi d\theta\\
=&\int_0^A\sin^2\Phi(\theta)d\theta\\
=&\int_0^A \tan^2b/(\tan^2b+\cos^2\theta)d\theta\\
=&\tan^2 b\int_0^A \sec^2\theta/(1+\tan^2b+\tan^2b\tan^2\theta)d\theta\\
=&\tan^2 b\int_0^{\tan A}1/(\sec^2b+u^2\tan^2b)du\\
=&\sin b\int_0^{\tan A\sin b}1/(1+u^2)du\\
=&\sin b(\tan^{-1}(\sin b\tan A))\\
=&a\sin b
\end{align*}

Where Claim \ref{parama} was used in the third line and the last line follows from Napier's rule \ref{Nap3}.\\
\end{proof}

{\rem Denote the area of $\triangle ABC$ as $[ABC]$. The centroid $G$ satisfies

\[[ABC]G=\int_{\triangle ABC}Xd\mu(X)\]

Let $d$ be the length of the altitude from angle $A$ to side $a$. Either $\triangle ABC$ is the union of two right triangles with common side $d$, or the difference. Either way, use the previous claim for the second line and Theorem \ref{SLOS} in the third of the following,

\begin{align*}2[ABC](G\cdot A)=&2M_A(\triangle ABC)\\=&a\sin d\\=&a\sin b\sin C\end{align*}

Since these quantities are independent of coordinates, the two other cases ($A\to B\to C$) follow similarly. If $A,B,C$ have no common great circle, $G$ is the unique point satisfying these equations. The formula

\begin{align}2[ABC]G\mathrm{det}(A,B,C)=\sum_{\mathrm{cyc}}(B\times C)a\sin b\sin C\label{CentroidFormula}\end{align}

is verified by taking dot products with each of $A,B,C$ separately. Recall the triple product $\mathrm{det}(A,B,C)=A\times B\cdot C$. Now, arbitrarily putting $A=(1,0,0)$ and $C_2=0$ and $\sigma(\triangle ABC)=1$, we see that $C\times B\cdot A=C_1B_2=\sin b(\sin c\sin A)$, the volume of the zonotope generated by $A,B,C$, is the ubiquitous quantity $n$ associated to a spherical triangle in \cite{IT}. Dividing the equation \ref{CentroidFormula} thru by this quantity becomes Brock's formula \cite{JB},

\[2[ABC]G=\sum_\mathrm{cyc} (B\times C)\frac{a}{\sin a}\]}

Consider a right spherical triangle $\triangle$ with $C=\pi/2$. We use equation \ref{Nap1} to convert the formula from Theorem \ref{avgheight3D} in terms of $A$ and $B$.

{\thm\label{NegDef3D} The function defined by \[f(A,B):=a\sin b=\left(\cos^{-1}\frac{\cos B}{\sin A}\right)\left(1-\frac{\cos^2A}{\sin^2B}\right)^{1/2}\] is negative definite (i.e. $-f$ is convex) in the region \[R:=\{-\pi/2<A,B<\pi/2\}\cap\{\cos^2A+\cos^2B<1\}\]}

\begin{proof} See section \ref{ProofNegDef3D}\end{proof}

{\thm The Simplex Mean Width Conjecture is true for $d=3$}\label{dim3thm}

\begin{proof}
Let the orientation of a spherical triangle be $\sigma(\triangle ABC)=+1$ for a counterclockwise ordering of $A,B,C$ and $-1$ for clockwise. If $\sigma(\triangle ABC)=-1$, measure both the angles and edges as negative.\\

Suppose $\triangle=\mathrm{conv}(v_0,\dots, v_3)$ is any simplex with $|v_i|=1$ for each $i$ (by Claim \ref{basic}). The Voronoi cells

\[S_i(\triangle):=\{u\in\sph:i\in\mathrm{argmax}_i u\cdot v_i\}\]

are each the intersection of three hemispheres and so are spherical triangles tiling $\sph$. For each $i$, label the vertices of $S_i$ as $A,B,C$ (suppressing subscript $i$) in the counterclockwise orientation. Draw arcs from $v_i$ to $A,B,C$ and drop altitudes from $v_i$ to the edges $a,b,c$ bounding $S_i$ (they may leave the interior of $S_i$) with feet $D,E,F$ on $a,b,c$ (resp.). Form the collection $\calT$ of right triangles $\triangle v_iXY$ such that $X,Y$ come from the list $A,F,B,D,C,E$ and $Y$ is immediately to the right of $X$. Since

\[\sigma(\triangle XYZ)\mathbbm{1}_{\triangle XYZ}+\sigma(\triangle XZW)\mathbbm{1}_{\triangle XZW}=\sigma(\triangle XYW)\mathbbm{1}_{\triangle XYW}\]

we can write

\[M_{v_i}(S_i)=\sum_{T\in\calT}\sigma(T)M_{v_i}(T)\]

Also, the signed sum of the six angles meeting at $v_i$ is $2\pi$ while the signed sum of the other six non-right angles equals the sum of the angles of $S_i$, which exceeds the area of $S_i$ by $\pi$ (Girard's Theorem). It follows that the whole of $\sph$ may be divided into 24 right triangles where the signed sum of the 24 angles measured at $\{v_0,\dots,v_3\}$ is $4*2\pi$ and the sum of the other 24 non-right angles is $4\pi+4\pi$. To see that none of these 24 right triangles contains an angle exceeding $\pi/2$, note that the spherical Pythagorean (\ref{SPT}) implies there would be two such edges from the same triangle, in particular one meeting at a vertex $v_i$. But, from (Claim. \ref{basic}(c)) every point in $\sph$ is within $\pi/2$ from a vertex $v_i$. As no edge exceeds $\pi/2$, neither does any angle. Theorems \ref{avgheight3D} and \ref{NegDef3D} imply that the maximum mean width occurs when all these angles are equal, i.e. a regular tetrahedron.
\end{proof}

\section{Integration on Spheres}\label{IntOnSph}
To generalize Theorem \ref{avgheight3D} we need to generalize a right triangle to a right angled simplex, called a path simplex, on higher dimensional spheres. Then, we find the marginal mean of an RV uniform over the simplex. Let us start with the measure of spherical caps. Denote a spherical cap of geodesic radius $r$, centered at $e_1$, by $B_\mathbb{S}^{d-1}(r)\subset\sphd$. Define the incomplete Wallis integral by

\[\mathcal{W}^d(x)=\int_0^x \sin^dtdt\]

and denote $W^d:=\mathcal{W}^d(\pi)$.

{\thm The uniform measure of a spherical cap is \begin{align}\mu(B_\mathbb{S}^{d-1}(r))=\frac{\calW^{d-2}(r)}{W^{d-2}}\label{CapMeas}\end{align} where \begin{align}\calW^d(r)=\frac{d-1}{d}\calW^{d-2}(r)-\frac{1}{d}\cos r\sin^{d-1} r\label{WallisRecursion}\end{align}}

\begin{proof}
Cut the sphere into thin slices perpendicular to the radius at the center of the cap. The cross section at geodesic radius $\varphi$ is $\sin\varphi\cdot\mathbb{S}^{d-2}$, so the measure of a thin slice is proportional to $\sin^{d-2} \varphi d\varphi$, and equation \ref{CapMeas} follows.\\

For the recursion formula, consider $(\sin^d \varphi)^{\prime\prime}$ for $d\ge 2$:

\begin{align*}(d\cos \varphi\sin^{d-1}\varphi)^\prime=&d(d-1)\cos^2 \varphi\sin^{d-2}\varphi-n\sin^d \varphi\\=&d(d-1)\sin^{d-2}\varphi-d^2\sin^d \varphi\end{align*}

Integrating gives \ref{WallisRecursion}.
\end{proof}
{\rem Equation \ref{WallisRecursion} may be used to find the trigonmetric series for $\calW^d(r)$. Note that $\calW^0(r)=r$ and $\calW^1(r)=1-\cos r$. The denominator of \ref{CapMeas} is the well known Wallis Integral. When $r=\pi$, we can multiply equation \ref{WallisRecursion} through by $W^{d-1}$ to get

\[d\cdot W^dW^{d-1}=(d-1)W^{d-1}W^{d-2}\]

which shows \[d\cdot W^dW^{d-1}=W^1W^0=2\pi\] for all $d$. Since $W^d$ is decreasing with $d$, we come upon the bound (see also \cite{ST}),

\begin{align}\sqrt{\frac{2\pi}{d+1}}<&W^d<\sqrt{\frac{2\pi}{d}}\end{align}}

{\cor Averaging the marginal of an RV uniform over a cap,

\begin{align*}M_{e_1}(B_{\mathbb{S}}^{d-1}(r))=&\int_0^r \cos t\ d\mu(B_{\mathbb{S}}^{d-1}(t))\\=&\frac{\int_0^r\cos t\sin^{d-2} tdt}{W^{d-2}}\\ =&\frac{\sin^{d-1} r}{(d-1)W^{d-2}}\end{align*}}

Now, fix an axis $A\in\sphd$ and parametrize $\sphd$ by $X=(r, \theta)\in [0,\pi]\times\mathbb{S}^{d-2}$ where $\cos r=X\cdot A$ and 

\begin{align*}\theta=\frac{X-A\cos r}{|X-A\cos r|}\label{Param}\end{align*}

{\cor Let a region $R\subset \sphd$ be defined by $0\le r\le \varphi(\theta)$ under the above parameterization for some measurable $\varphi$. Then

\begin{align*}M_A(R)=\int_{\mathbb{S}^{d-2}} \frac{\sin^{d-1}\varphi(\theta)}{(d-1)W^{d-2}} d\mu(\theta)\end{align*}}\label{MAR}

When this region is a spherical simplex, say $T\subset\sphd$, it is the intersection of $d$ hemispheres. Arrange the vertices as columns in the matrix $V$. We describe $T$ with matrix $H\in \RR^{d\times d}$, whose rows are the centers of the $d$ hemispheres, ordered so that the $ith$ vertex is in the interior of the $i$th hemisphere. Then $x/|x|\in T$ iff

\[Hx\ge\mathbf{0}\]

Since each vertex $v_i$ is on the boundary of the other $d-1$ hemispheres, $H$ is a diagonal multiple of $V^{-1}$. Further, if $e_1$ is a vertex, we may assume the first column of $H$ is $(1,0,0,\dots, 0)^t$. Also, $T$ is a union of arcs originating at $e_1$, and $T$ has a parameterization $0\le r\le \varphi(\theta)$ as in Corollary \ref{MAR}. To find $\varphi(\theta)$, note that the first row of $H$ describes the facet opposite $e_1$. 

\[\cos \varphi(\theta)+\theta\cdot(H_{12},\dots,H_{1d})\sin \varphi(\theta) = 0\]

giving

\[\tan\varphi(\theta)=-\frac{1}{\theta\cdot(H_{12},\dots,H_{1d})}\]

or

\[\sin^2\varphi(\theta)=\frac{1}{1+(\theta\cdot(H_{12},\dots,H_{1d}))^2}\]

The support ${\tilde T}$ of $\varphi$ is the simplex described by omitting the first row and column of $H$. Applying Corollary \ref{MAR}, we get

\begin{align}M_{e_1}(T)=\frac{1}{(d-2)W^{d-2}}\int_{{\tilde T}}(1+(\theta\cdot(H_{12},\dots,H_{1d}))^2)^{(1-d)/2}d\mu(\theta)\label{MAT}\end{align}

\section{SMWC For $d>3$}\label{dimd}

In this section, we assume the set of points $V=\{v_0,\dots,v_d\}\subset\sphd$ is in general position so that the feet of the altitudes to a face never land in the boundary of that face.

{\defin A \emph{path simplex} in $\sphd$ is a spherical simplex with $d$ edges forming a path such that any pair of them determine great circles meeting at right angles. The endpoints of the path are called \emph{end vertices}.}\label{DefPathSimplex}

{\defin Let $T$ be a spherical simplex with vertices $\{p_1,\dots,p_d\}$ and inward facing normals $\{g_1,\dots,g_d\}$ corresponding to the opposite facets. Let $G_T=[g_1\dots g_d]$ and the \emph{(angle) Gram matrix} be $G_T^tG_T$.}\label{GramDef}\\

Note, many authors use $2I_{d\times d}-G_T$ as the Gram matrix. It is known that $G$ determines $T$ up to isometry, e.g. \cite{RD}.

{\claim If $T$ is a path simplex in $\sphd$, it has a tridiagonal Gram matrix.}
\begin{proof}
Say the vertices along the path are $p_1,\dots, p_d$, in order. The $i$th row of $P^{-1}=[p_1,\dots,p_d]^{-1}$ is normal to $F_i$, so there is a diagonal matrix $D$ with

\begin{align}G=(DP^{-1})(DP^{-1})^t=D(P^tP)^{-1}D\label{GramFromV}\end{align}

The path simplex condition may be restated as the arc $p_ip_{i+1}$ must be orthogonal to $\mathrm{span}\{p_{i+1},\dots,p_d\}$ for $1\le i\le d-1$. That is, the vector $p_i-(p_i\cdot p_{i+1})p_{i+1}$ tangent to the arc at $p_{i+1}$ is orthogonal to each of $p_{i+1},\dots, p_d$. It follows that

\[P^tP\left(\begin{array}{ccccc}1 & 0 & 0 & \dots & 0\\ p_1\cdot p_2 & 1 & 0 & \dots & 0\\ 0 & p_2\cdot p_3 & 1 & \dots & 0\\ \vdots & \vdots & \vdots & \ddots & \vdots\\ 0 & 0 & 0 & \dots & 1\end{array}\right)\]

is upper triangular. Substitute this into equation \ref{GramFromV}. Since $G$ is symmetric, $G$ must be tridiagonal.
\end{proof}

{\claim Let $\triangle$ be spherical simplex in $\sphd$. For almost every $O\in\sphd$, there is a collection $\calT$ of $d!$ path simplices with end vertex $O$ and a sign $\sigma:\calT\to\{-1,1\}$ satisfying

\begin{align}\mathbbm{1}_\triangle=\sum_{T\in\calT}\sigma(T)\mathbbm{1}_T\hspace{.5in}(a.e.)\end{align}}\label{IndicatorPartition}

\begin{proof}
Proceed by induction; $d=2$ is trivial. Suppose the lemma is true at level $d$, and we will show it is also true at the level $d+1$.\\

Let $\triangle\subset \mathbb{S}^{d}$ be a simplex. Drop an altitude from $O\in\sphd$ to any facet $F$ of $\triangle$ with foot $O_F\not\in\partial F$ (excluding a measure $0$ case). Both $F$ and $O_F$ lie in a common great sphere $S_F$ of dimension $d-1$. By induction, we have a collection $\calT_F$ of $d!$ path simplices with paths beginning at $O_F$ and satisfying

\[\mathbbm{1}_{F}=\sum_{T_F\in\calT_F}\sigma(T_F)\mathbbm{1}_{T_F}\]

Since $OO_F$ is perpendicular to $S_F$, each of these paths may append $OO_F$ to the beginning to represent the path simplices $T^F:=\mathrm{conv}(O,T_F)$. Form the collection $\calT^F:=\{T^F:T_F\in\calT_F\}$, and the simplex $\triangle^F:=\mathrm{conv}(O,F)$. Now, identify points $X\in S_F$ with segments $OX$, so that

\begin{align}
\mathbbm{1}_{\triangle^F}=&\sum_{T^F\in \calT^F}\sigma( T_F)\mathbbm{1}_{T^F}\nonumber\\
\mathbbm{1}_{\triangle}=&\sum_{F}\left.\begin{cases}-1 &\text{ if } \mu(\triangle^F\cap \triangle)=0\\ +1&\text{ o.w.}\end{cases}\right\}\mathbbm{1}_{\triangle^F}
\end{align}

Since there are $d+1$ faces $F$ (in the second line above), the total sum involves $(d+1)!$ path simplices.
\end{proof}

{\claim The Voronoi cells generated by $V$ are simplices. The set $\mathcal{C}$ of faces across all $d+1$ cells corresponds to the nonempty subsets of $V$ in that each face is the set of points equidistant from the points in its corresponding subset. Also, the intersection of any two of these faces is again a face ($\mathcal{C}$ is a simplicial complex).}\label{SubsetCorr}

\begin{proof}
The Voronoi cells are bounded by perpendicular bisectors between each pair of vertices $v_i,v_j$. That's $d$ facets for each cell, i.e. a simplex. Each face of a simplex corresponds to an intersection of facets, or in this case, perpendicular bisectors. Thus, each facet is the set of points equidistant from some subset of $V$. It follows that the intersection of two faces corresponds to the union of the two corresponding subsets.
\end{proof}

For a set $V=\{v_0,\dots, v_d\}\subset \sphd$ and a maximal chain $\tau$ in the power set of $V$, denote by $S_\tau$ the path simplex starting at the vertex in the singleton set of $\tau$ and landing successively in the faces determined by the subsets in $\tau$, increasing by inclusion.

{\claim Let $\mathcal{C}$ be as in the previous claim. To each face $F\in \mathcal{C}$ (and hence every subset $Q\subset V$) call $p(F)\in\mathrm{span}(F)$ the point that minimizes the distance from the points in $Q$. Then, $p(F)$ is a vertex of $S_\tau$ iff $Q\in \tau$.}\\

\begin{proof}
$\mathrm{span}(F)$ is the set of points equidistant from each point in $Q$. Since the composition of projections onto subspaces is the projection onto the intersection of the subspaces, each point in $Q$ projects onto $\mathrm{span}(F)$ with foot $p(F)$. Finally, $p(F)$ is a vertex of $S_\tau$ if $p(F)$ is equidistant from points in $Q$, but no other points in $V$.
\end{proof}

The following is the main result.\\

{\thm Let $T\subset\sphd$ be the path simplex with $\theta_{ij}$ the angle between facets $F_i$ and $F_j$. If $A$ is an end vertex of $T$ and the Hessian of $M_A(T)$ from equation \ref{MAT} as a function of $\theta_{12},\dots,\theta_{d,d+1}$ is negative definite, then the simplex mean width conjecture is true.}

\begin{proof}
Let $\triangle$ be a simplex in $\mathbb{R}^d$ with vertices $V=\{v_0,\dots,v_d\}\subset\sphd$ (Claim \ref{basic}(b)). The altitudes from $v_i$ and $v_j$ to the perpendicular bisector between them meet at the same point. The rest of the triangulations of each cell, as in Claim 4.4, coincide along facets, forming a $d+1$ times larger simplicial complex $\mathcal{C}^\dagger$ where every face is the intersection of cells.\\

By Claim 4.6 the vertices (or facets or acute dihedral angles) of any path simplex $S_\tau$ may be ordered as $\tau$ is ordered. Given three consecutive subsets $Q,Q\cup\{x\}, Q\cup\{x,y\}$ in a maximal chain $\tau$ of the subsets of $V$, the only replacement subset for $Q\cup\{x\}$ to keep $\tau$ a chain is $Q\cup\{y\}$. It follows that if a dihedral angle at a $d-3$ dimensional face in a path simplex $S_\tau$ is between two nonadjacent facets, there are $2\times 2=4$ path simplices meeting at that face, at right angles as we know from Claim 4.3. Similarly, a dihedral angle between two adjacent facets has $3\times 2=6$ path simplices meeting at that face.\\

$\mathcal{C}^\dagger$ has a total of $(d+1)!$ path simplices, each with $d$ acute dihedral angles meeting together at the $d-3$ dimensional faces, according to their position in the $\tau$ ordering. So, for each $i$, the sum of the angles $\theta_{i,i+1}$ over all path simplices is $2\pi (d+1)!/6$. If $M_A(T)$ is negative definite as a function of these adjacent dihedral angles, the mean width from equation \ref{SumMAR} will be maximum when all the adjacent dihedral angles are equal to $\pi/3$, that is when $\triangle$ is regular.
\end{proof}

\section{Proof of Theorem 2.3}\label{ProofNegDef3D}

 Use Napier's formula \ref{Nap1} to phrase quantities in terms of `intermediate variables' $a,b$. Differentiating this formula with respect to $A$ and using the spherical law of sines \ref{SLOS},

\begin{align*}
\sin A=\sin a\sin B\frac{da}{dA}=\sin A\sin b\frac{da}{dA}
\end{align*}

Now with respect to $B$, and using again \ref{SLOS} and \ref{Nap1} (cycled $A\to B$),

\begin{align*}
0=&\sin a\sin B\frac{da}{dB}-\cos a\cos B\\
=&\sin A\sin b\frac{da}{dB}-\cos a\cos b\sin A
\end{align*}

The same holds if $(A,a)$ and $(B,b)$ are interchanged, yielding

\begin{align*}
  \begin{split}
    \frac{\partial a}{\partial A}=&\ \frac{1}{\sin b}\\
    \frac{\partial b}{\partial A}=&\ \frac{\cos a\cos b}{\sin a}\\
  \end{split}
  \begin{split}
    \frac{\partial a}{\partial B}=&\ \frac{\cos a\cos b}{\sin b}\\
    \frac{\partial b}{\partial B}=&\ \frac{1}{\sin a}
  \end{split}
\end{align*}

The first order partials are then

\begin{align*}
f_A=&\frac{da}{dA}\sin b+a\cos b\frac{db}{dA}=1+\frac{a\cos a\cos^2 b}{\sin a}\\
f_B=&\frac{da}{dB}\sin b+a\cos b\frac{db}{dB}=\cos b\left(\cos a+\frac{a}{\sin a}\right)
\end{align*}

For second order partials,

\begin{align*}
\frac{\sin a}{\cos a\cos^2 b}f_{AA}=&\frac{1}{\sin b}-\left(2\sin b\cos a+\frac{1}{\sin b\cos a}\right)\frac{a}{\sin a}\\
\frac{\sin a}{\cos a}f_{BB}=&-\sin b+\frac{\cos^2a\cos^2b}{\sin b}-\left(\frac{\sin b}{\cos a}+\frac{\cos^2 b}{\sin b}\cos a\right)\frac{a}{\sin a}\\
\frac{\sin a}{\cos a\cos b}f_{AB}=&\frac{\cos^2 b}{\sin b}\cos a-\left(\frac{1+\sin^2 b}{\sin b}\right)\frac{a}{\sin a}
\end{align*}

Uniquely represent

\begin{align*}
\frac{\sin^2 a}{\cos^2a\cos^2b}&\mathrm{detHess}(f)=\frac{\sin^2 a}{\cos^2a\cos^2b}(f_{AA}f_{BB}-f_{AB}^2)\\
=&P(\sin a,\sin b)+Q(\sin a,\sin b)\frac{a\cos a}{\sin a}+R(\sin a,\sin b)\frac{a^2}{\sin^2 a}
\end{align*}

for some rational functions $P,Q,R$. Now we will compute P,Q,R.

\begin{align*}
P(x,y)=&\left(\frac{1}{y}\right)\left(-y+\frac{(1-x^2)(1-y^2)}{y}\right)-\left(\frac{(1-y^2)^2}{y^2}(1-x^2)\right)\\
=&-1+\frac{(1-x^2)(1-y^2)(1-(1-y^2))}{y^2}\\
=&(1-x^2)(1-y^2)-1\\
Q(x,y)=&-\frac{1}{y}\left(\frac{y}{1-x^2}+\frac{1-y^2}{y}\right)\\
&-\left(2y+\frac{1}{y(1-x^2)}\right)\left(-y+\frac{(1-x^2)(1-y^2)}{y}\right)\\
&+2\left(\frac{1-y^2}{y}\right)\left(\frac{1+y^2}{y}\right)\\
=&2\left(y^2-\frac{1-y^2}{y^2}-(1-x^2)(1-y^2)+\frac{1-y^4}{y^2}\right)\\
=&2(1-(1-x^2)(1-y^2))\\
R(x,y)=&\left(2y+\frac{1}{y(1-x^2)}\right)\left(y+\frac{(1-y^2)(1-x^2)}{y}\right)-\left(y+\frac{1}{y}\right)^2\\
=&2y^2+\frac{1}{1-x^2}+2(1-x^2)(1-y^2)+\frac{1}{y^2}-1-\left(y+\frac{1}{y}\right)^2\\
=&-y^2-1-2x^2+2x^2y^2+\frac{1}{1-x^2}\\
=&-2+(1-2x^2)(1-y^2)+\frac{1}{1-x^2}\\
=&(2x^2-1)\left(y^2-1+\frac{1}{1-x^2}\right)\\
=&\left(\frac{x^2}{1-x^2}-\frac{1-x^2}{x^2}\right) x^2\left(1-(1-x^2)(1-y^2)\right)
\end{align*}

From the spherical law of sines \ref{SLOS} and spherical Pythagorean \ref{SPT},

\[\frac{\sin^2 a}{\sin^2 A}=\sin^2 c=1-\cos^2 c=1-\cos^2a\cos^2 b\]

and so we reduce the Hessian to

\begin{align*}
\frac{\sin^2 A}{\cos^2 a\cos^2 b}\mathrm{detHess}f=&(\tan^2 a-\cot^2 a)a^2 +(2\cot a)a-1\\
=&a^2\tan^2 a-(1-a\cot a)^2
\end{align*}

This can be seen to be positive by following the chain of inequalities backward:

\begin{align*}
0>&(1-a\cot a-a\tan a)(1-a\cot a+a\tan a)\\
=&\left(1-\frac{a}{\sin a\cos a}\right)\left(1-a\frac{\cos^2 a-\sin^2a}{\sin a\cos a}\right)\\
=&\left(1-\frac{2a}{\sin 2a}\right)\left(1-\frac{2a}{\tan 2a}\right)
\end{align*}

for $|a|<\pi/2$, which follows since $x/\sin x$ (resp. $x/\tan x$) is tangent to $y=1$ and increasing (resp. decreasing) away from $0$ for $-\pi<x<\pi$.

Thanks to Joshua Horowitz and the use of his copy of Mathematica, which cleared the way for this proof.

\end{document}